\documentclass[12pt,leqno]{amsart}
\usepackage{amsfonts,amsmath,amssymb,amsthm, geometry}

\usepackage[hypertex]{hyperref}
\usepackage[all]{xy}

\newtheorem{theorem}{Theorem}[section]
\newtheorem{corollary}[theorem]{Corollary}
\newtheorem{lemma}[theorem]{Lemma}

\theoremstyle{definition}
\newtheorem{definition}[theorem]{Definition}
\theoremstyle{remark}
\newtheorem{remark}[theorem]{\sc Remark}
\newtheorem{example}[theorem]{\sc Example}
\begin{document}
\title[Local fibrations of real map germs]{Equivalence of real Milnor's fibrations for quasi homogeneous singularities}

\author[A. dos Santos]{Ara\'ujo dos Santos, R.}

\email{rnonato@icmc.usp.br}
\address[R. N. Ara\'ujo dos Santos]{Universidade de S\~ao Paulo, Instituto de Ci\^encias Matem\'aticas e de Computa\c c\~ao
, Av. Trabalhador S\~ao-Carlense, 400 -
Centro, Postal Box: 668, S\~ao Carlos - S\~ao Paulo - Brazil, Postal Code: 13560-970.}

%\author{Mihai}
%\address[Mihai]{somewhere}
\date{\today}
\maketitle
%%%%%%%%%%%%%%%%%%%
\begin{abstract}

We are going to use the Euler's vector fields in order to show that for real quasi-homogeneous singularities with isolated critical value, the Milnor's fibration in a "thin" hollowed tube involving the zero level and the fibration in the complement of "link" in sphere are equivalents, since they exist. Moreover, in order to do that, we explicitly characterize the critical points of projection $\frac{f}{\|f\|}:S_{\epsilon}^{m}\setminus K_{\epsilon}\to S^{1}$, where $K_{\epsilon}$ is the link of singularity.

\end{abstract}

\section{Introduction}

It is well known (see for instance \cite{Mi}) the Milnor studies about fibrations in a neighborhood of a (isolated)critical point. For real isolated mappings, he showed that if $f:\mathbb{R}^{m},0\to \mathbb{R}^{p},0, m\geq p\geq 2$, is a real polynomial map germ whose derivative $D{f}$
has rank $p$ on a punctured neighborhood of $0\in \mathbb R^m$,
then there exists $\epsilon >0$ and $\eta >0$, sufficiently small
$0<\eta \ll\epsilon <1$, such that considering
$E:=\overline{B_{\epsilon}^{m}(0)}\cap f^{-1}( S_{\eta}^{p-1})$, where
$B_{\epsilon}^{m}(0)$ is the open ball centered in $0\in
\mathbb{R}^{m}$ with radius $\epsilon $, we have that

\begin{center}

$\displaystyle{f_{\mid_{E}}: E \to S_{\eta}^{p-1}}$ {\bf (1)}

\end{center}

is a smooth locally trivial fibre bundle. Furthermore, using the flow of a vector fields $\omega(x)$ in $\overline{B_{\epsilon}^{m}(0)}\setminus \{f=0\}$ satisfying the conditions:

\begin{equation*}
\begin{cases}
\langle \omega(x), x \rangle >0; \\
\langle \omega(x), \nabla (\|f(x)\|^{2}) \rangle >0,
\end{cases}
\end{equation*}

Milnor constructed a diffeomorphism
that ``pushes" $E$ to $S^{m-1}_{\epsilon}\setminus N_{K_{\epsilon}}$,
where $N_{K_{\epsilon}}$ denotes an open tubular neighborhood of the
link $K_{\epsilon}$ in $S^{m-1}_{\epsilon}.$ Moreover, this fibre
bundle can be extended to the complement of the link $K_{\epsilon}$ in the sphere
$S^{m-1}_{\epsilon}$, providing this way a fibration

\begin{center}

$S^{m-1}_{\epsilon}\setminus K_{\epsilon}\to S^{p-1}$. {\bf (2)}

\end{center}

As is obviously seen in example below,
in all these constructions we cannot guarantee that the natural map
$\displaystyle{\frac{f}{\|f \|}}$ is the projection
of the fibration(see details in example.\ref{ex2}):

\begin{example}\cite[page 99]{Mi}
\begin{equation*}
\begin{cases}
P=x \\
Q=x^{2}+y(x^{2}+y^{2})
\end{cases}
\end{equation*}
\end{example}

In what follow we are going to characterize the critical points of the projection map $\displaystyle{\frac{f}{\|f\|}:{S}^{m-1}_{\epsilon}\setminus K_{\epsilon}\to {S}^{1}}$ (see Lemma.\ref{p1}) in order to understand whether the fibrations {\bf (1)} and {\bf (2)} are equivalents(see definition \ref{de}). In fact, we will prove that for real non-isolated quasi-homogeneous singularities these two fibrations are equivalents, since they exists.

\section{Definitions and setup}

\begin{definition}\cite{RSV}

Let $\psi:(\mathbb{R}^{m},0)\to (\mathbb{R}^{p},0)$, $m\geq p\geq
2$, be a map germ with isolated singularity at the origin. If for
all $\epsilon >0$ sufficiently small, the map $\phi
=\displaystyle{\frac{\psi}{\|\psi\|}}:S^{m-1}_{\epsilon}\setminus
K_{\epsilon} \,\to \,{S}^{p-1}\,$, is a projection of a smooth
locally trivial fibre bundle, where $K_{\epsilon}$ is the link of
singularity at $0$(or empty), we say that the map germ  satisfies the {\bf
Strong Milnor condition} at origin $0\in \mathbb{R}^{m}.$
\end{definition}

This problem first appeared in \cite{Ja1}, \cite{Ja2}, and further have been developed in \cite{RSV}, \cite{Se1}, \cite{Se2}, \cite{RS}, \cite{SM}.

\vspace{0.5cm}

\noindent Consider $f:\mathbb{R}^{m},0\to \mathbb{R}^{2},0 $, $f(x)=(P(x),Q(x))$, where $P,Q:\mathbb{R}^{m},0\to \mathbb{R},0$ are polynomial functions. Define $\omega (x)=P(x)\nabla Q(x)-Q(x)\nabla P(x)$ and ${\bf e}(x)=({w_{1}}x_{1},..., {w_{m}}x_{m})$, $w_{i}\in \mathbb{Q}^{*}_{+}$ for all $i=1,...,m$. Suppose that $P,Q:\mathbb{R}^{m}\to \mathbb{R}$ are quasi(weighted-) homogeneous of type $(w_{1},...,w_{m}; b_{1})$ and $(w_{1},...,w_{m}; b_{2})$, i.e, $P(\lambda^{w_{1}}x_{1},..., \lambda^{w_{m}}x_{m})=b_{1}P(x)$ and $Q(\lambda^{w_{1}}x_{1},..., \lambda^{w_{m}}x_{m})=b_{2}Q(x)$, for some $\lambda \in \mathbb{R}^{*}_{+}$. In this case we say that the functions $P$ and $Q$ are quasi-homogeneous of same weight and degrees $b_{1}$ and $b_{2}$, respectively. The vector fields ${\bf e}(x)$ is called the Euler's vector fields to $P$ and $Q$. This well known vector fields satisfies the following property:

\vspace{0.5cm}

\begin{center}

$\langle \nabla P(x), {\bf e}(x)\rangle =b_{1}P(x)$, $\langle \nabla Q(x), {\bf e}(x)\rangle = b_{2}Q(x)$.

\end{center}

\vspace{0.5cm}

In this paper I will show that for quasi-homogeneous map $f:\mathbb{R}^{m},0\to \mathbb{R}^{2},0$, of same weight and degree, the two fibrations above, {\bf (1)} and {\bf (2)} are equivalents. Moreover, the diffeomorphism which provide this equivalence is given by the flow of Euler's vector fields.

The main result is:

\begin{theorem}\label{HR}

Let $f=(P,Q):\mathbb{R}^{m},0\to \mathbb{R}^{2},0$ quasi-homogeneous polynomial map germ of same type and degree, with an isolated critical value at origin $0\in \mathbb{R}^{2}$. Suppose that the two fibrations, {\bf (1)} and {\bf (2)} exists. Then, they are equivalents.

\end{theorem}

\begin{corollary}
Consider $h:\mathbb{C}^{m},0\to \mathbb{C},0$, $h(z)=f(z)\overline{g(z)}$, where $f,g:\mathbb{C}^{m}\to \mathbb{C}$ are holomorphic functions. Suppose that $0\in \mathbb{C}$ is a isolated critical value of $h$ and the pair $(Re(h(z)), Im(h(z)))$ are quasi-homogeneous of same weight and degree. Then, the two fibrations {\bf (1)} and {\bf (2)} are equivalents.

\end{corollary}

This kind of fibration have been further studied in \cite{Se2} and \cite{pi}.

\section{Critical points of projection}

In what follows we are considering a real analytic map germ $\displaystyle{f=(P,Q): (\mathbb{R}^{m},0)\to (\mathbb{R}^{2},0)}$, and the projection $\displaystyle{\frac{f}{\|f\|}:{S}^{m-1}_{\epsilon}\setminus K_{\epsilon}\to {S}^{1}}$.

\vspace{0.3cm}

\noindent Notation:

\noindent For a general smooth function $g:\mathbb{R}^{m}\to \mathbb{R}$, for convenience we will denote the derivative $\displaystyle{\frac{\partial g(x)}{\partial x_{i}}}$ by $\partial_{i}g(x)$;

\vspace{0.3cm}

\noindent In the sequence we shall see that the vector fields $\omega (x)=P(x)\nabla Q(x)-Q(x)\nabla P(x)$above shows a important role in order to describe the critical points of the map projection $\displaystyle{\frac{f}{\|f\|}}$. The lemma below appeared first in \cite{Ja2}, but here we are given a different proof.

\begin{lemma}\cite{Ja2}\label{p1}

A point $x\in {S}^{m-1}_{\epsilon }\setminus K_{\epsilon }$ is a critical point of $\displaystyle{\frac{f}{\|f\|}}$ if, and only if, the vector $\displaystyle{\omega(x) =\lambda x}$, for some real $\lambda $.

\end{lemma}

\proof We know that $\displaystyle{\frac{f(x)}{\|f(x)\|}=(\frac{P(x)}{\sqrt{P(x)^{2}+Q(x)^{2}}},\frac{Q(x)}{\sqrt{P(x)^{2}+Q(x)^{2}}} )}$ and the jacobian matrix is

\vspace{0.2cm}

$\displaystyle{J(\frac{f(x)}{\|f(x)\|})=\left(
 \begin{array}{c}
 \nabla (\frac{P(x)}{\sqrt{P(x)^2+Q(x)^2}}) \\
 \nabla (\frac{Q(x)}{\sqrt{P(x)^2+Q(x)^2}}) \\
 \end{array}
\right) }$.

\vspace{0.2cm}

So, using notation above we have:

\vspace{0.3cm}

$\displaystyle{\partial_{i} (\frac{P(x)}{\sqrt{P(x)^{2}+Q(x)^{2}}}
)=\frac{\partial_{i}(P(x))\sqrt{P(x)^{2}+Q(x)^{2}}-P(x)\partial_{i}(\sqrt{P(x)^{2}+Q(x)^{2}})}{P(x)^{2}+Q(x)^{2}} }$ {\bf $(*)$}.

\vspace{0.2cm}

The second derivative on numerator is $\displaystyle{\frac{2(P(x)\partial_{i}P(x)+Q(x)\partial_{i}Q(x))}{2\sqrt{P(x)^{2}+Q(x)^{2}}}}$. {\bf $(\#)$}

Combining {\bf $(*)$} and {\bf $(\#)$} we have:

\vspace{0.3cm}

$\displaystyle{\partial_{i}(\frac{P(x)}{\sqrt{P(x)^{2}+Q(x)^{2}}}
)= \frac{\partial_{i}(P(x))(P(x)^{2}+Q(x)^{2})-P(x)((P(x)\partial_{i}P(x)+Q(x)\partial_{i}Q(x)))}{(P(x)^{2}+Q(x)^{2})^{\frac{3}{2}}}}=$

$\displaystyle{\frac{Q(x)^{2}\partial_{i}P(x)-P(x)Q(x)\partial_{i}Q(x)}{(P(x)^{2}+Q(x)^{2})^{\frac{3}{2}}}}=\displaystyle{\frac{-Q(x)}{(P(x)^{2}+Q(x)^{2})^{\frac{3}{2}}}(P(x)\partial_{i} Q(x)-Q(x)\partial_{i} P(x))},$ for all $i=1,...,n$

It means that the first line of jacobian matrix is

$\displaystyle{\frac{-Q(x)}{(P(x)^{2}+Q(x)^{2})^{\frac{3}{2}}}(P(x)\nabla Q(x)-Q(x)\nabla P(x))}$,
or $\displaystyle{\frac{-Q(x)}{(P(x)^{2}+Q(x)^{2})^{\frac{3}{2}}}\omega (x)}$.

\vspace{0.3cm}

Similarly, we have

$\displaystyle{\partial_{i} (\frac{Q(x)}{\sqrt{P(x)^{2}+Q(x)^{2}}}
)=\frac{\partial_{i}(Q(x))\sqrt{P(x)^{2}+Q(x)^{2}}-Q(x)\partial_{i}
(\sqrt{P(x)^{2}+Q(x)^{2}})}{P(x)^{2}+Q(x)^{2}} }=$

\vspace{0.3cm}

$=\displaystyle{\frac{\partial_{i}(Q(x))(P(x)^{2}+Q(x)^{2})-Q(x)((P(x)\partial_{i}P(x)+Q(x)\partial_{i}Q(x)))}{(P(x)^{2}+Q(x)^{2})^{\frac{3}{2}}}}=$

\vspace{0.3cm}

$=\displaystyle{\frac{P(x)^{2}\partial_{i}Q(x)-Q(x)P(x)\partial_{i}P(x)}{(P(x)^{2}+Q(x)^{2})^{\frac{3}{2}}}=(\frac{P(x)}{(P(x)^{2}+Q(x)^{2})^{\frac{3}{2}}})\omega (x)}$.

\vspace{0.3cm}

It means that the jacobian matrix have the following lines

$\displaystyle{J(\frac{f(x)}{\|f(x)\|})=\left(
\begin{array}{c}
 (\frac{-Q(x)}{(P(x)^{2}+Q(x)^{2})^{\frac{3}{2}}})\omega (x) \\
 (\frac{P(x)}{(P(x)^{2}+Q(x)^{2})^{\frac{3}{2}}})\omega (x) \\
 \end{array}
\right)}.$

\vspace{0.2cm}

Of course that, $\displaystyle{J(\frac{f(x)}{\|f(x)\|}): T_{x}({S}^{m-1}\setminus K_{\epsilon})\to T_{f(x)}{S}^{1}\simeq \mathbb{R}}$, so $x\in {S}^{m-1}\setminus K_{\epsilon}$ is a critical point of $\displaystyle{\frac{f(x)}{\|f(x)\|}}$ if, and only if, for all vectors $v\in T_{x}({S}^{m-1}\setminus K_{\epsilon})$ we have

\vspace{0.3cm}

$\displaystyle{J(\frac{f(x)}{\|f(x)\|}).v=
\left(\begin{array}{c}
  (\frac{-Q(x)}{(P(x)^{2}+Q(x)^{2})^{\frac{3}{2}}})\omega (x) \\
  (\frac{P(x)}{(P(x)^{2}+Q(x)^{2})^{\frac{3}{2}}})\omega (x) \\
   \end{array}
   \right).v=
\left(\begin{array}{c}
  (\frac{-Q(x)}{(P(x)^{2}+Q(x)^{2})^{\frac{3}{2}}})(\omega (x).v) \\
  (\frac{P(x)}{(P(x)^{2}+Q(x)^{2})^{\frac{3}{2}}})(\omega (x).v) \\
   \end{array}
\right)=
\left(
  \begin{array}{c}
    0 \\
    0 \\
  \end{array}
\right)}.$

\vspace{0.3cm}

Now, it is pretty easy to see that the critical locus of the map $\displaystyle{\frac{f(x)}{\|f(x)\|}}$ is precisely the points $x$ where $\omega (x)$ is parallel to vector position $x$. \endproof

\begin{example}\label{ex2}
\begin{equation*}
\begin{cases}
P=x \\
Q=x^{2}+y(x^{2}+y^{2})
\end{cases}
\end{equation*}

It is easy to see that $\nabla P(x,y)=(1,0)$, $\nabla Q(x,y)=(2x+2xy,x^{2}+3y^{2})$. So, $\omega (x,y)=(x^2+yx^2-y^3,x^3+3xy^2)$. The points $(x,y)$ where the vector fields $\omega (x,y)$ is parallel to vector position $(x,y)$ is equivalent to $\langle \omega (x,y), (-y,x)\rangle =0$. From the last equation we have $0=\langle \omega (x,y), (-y,x)\rangle =(x^2+y^2)^2-yx^2$. Now using polar coordinates, $x=r\cos(\theta), y=r\sin(\theta)$, in the last equation we have $r=\sin(\theta ).\cos^{2}(\theta )$, which is a non-degenerate curve going throughout origin.

\end{example}

\vspace{0.3cm}

\section{Preliminar results}

In order to proof main theorem we will do some small results which provide us extras informations. Let's see next that in our setting the projection is regular, i.e, have maximal rank.

\vspace{0.2cm}

\begin{lemma}\label{L1} Under conditions of Theorem \ref{HR} above we have $\langle \omega(x), {\bf e}(x) \rangle =0$.

\end{lemma}

\proof: It follows from easy calculation since\\

\noindent $\langle \omega (x), {\bf e}(x) \rangle = P(x)\langle \nabla Q(x), {\bf e}(x) \rangle - Q(x)\langle \nabla P(x), {\bf e}(x) \rangle = P(x)b Q(x)-Q(x)b P(x)=0$.\endproof

\begin{remark}\label{R1}

\noindent Since $\langle {\bf e}(x),x \rangle =\sum_{i=1}^{m}{w_{i}}x_{i}^{2}$, then $\langle {\bf e}(x),x \rangle =0$ iff $x=0.$

\end{remark}

\begin{lemma}\label{L2} Consider $V=\{ x\in {R}^{m}: P(x)=Q(x)=0\}=f^{-1}(0)$. Then, for all $x\in B_{\epsilon}(0)\setminus V$ the vector fields $\omega (x)$ and the vector position $x$ are not parallels.

\end{lemma}

 \proof If $\omega (x)$ were parallel to $x$ applying the Lemma \ref{L1} above we have $\langle x, {\bf e}(x) \rangle =0$. By remark \ref{R1}, this contradict the hypothesis that $x\in B_{\epsilon}(0)\setminus V$. \endproof

\vspace{0.3cm}

\begin{lemma} For $\epsilon $ small enough, there exist a smooth vector fields $\alpha $ on $B_{\epsilon}(0)\setminus V$ satisfying the following conditions:

\begin{enumerate}

\item $\langle \alpha (x), x\rangle >0;$

\item $\langle \alpha (x), \nabla (\|f(x)\|^{2}\rangle >0,$ where $\nabla (\|f(x)\|^{2})=2.(P(x)\nabla P(x)+ Q(x)\nabla Q(x)) $.

\item $\langle \alpha (x), \gamma (x)\rangle =0 $, \text{where} $\gamma (x)=P(x)\nabla Q(x)-Q(x)\nabla P(x)$.

\end{enumerate}

\end{lemma}

\proof Consider $\alpha (x):={\bf e}(x)=({w_{1}}x_{1},...,{w_{m}}x_{m})$ the Euler vector fields.

The conditions a) and c) follows from Lemma \ref{L1} and Remark \ref{R1}.

Condition b): Since $\nabla (\|f(x)\|^{2})(x)=2.(P(x)\nabla P(x)+Q(x) \nabla Q(x))$, then,

$\langle \nabla (\|(P,Q)\|^{2})(x),\alpha (x)\rangle = 2.(P(x)\langle \nabla P(x), \alpha (x)\rangle+Q(x) \langle \nabla Q(x), \alpha (x)\rangle)=2b.(P^{2}(x)+Q^{2}(x))$, where $b$ is the degree of quasi-homogeneity. Therefore, $\langle \nabla (\|(P,Q)\|^{2})(x),\alpha (x)\rangle >0$

\endproof

\section{Equivalence of Milnor fibrations}

\begin{definition}\label{de}

We say that two smooth locally trivial fibre bundle $f:M\to B$ and $g:N\to B$ are $\mathcal{C}^{k}-$equivalents(or $\mathcal{C}^{k}-$isomorphic), $k=0,1,...,\infty$ if there exist a $\mathcal{C}^{k}-$diffeomorphism $h$ such that the following diagram is commutative:

%%%%%%%%%%%%%%%%%%%%%%%%%%%%%%%%% commutative diagram%%%%%%%%%%%%%%%%%%

\begin{center}
$\xymatrix{\displaystyle{M}
\ar[d]_{f} \ar[r]^{h}
& \displaystyle{N} \ar[ld]^{g} \\
\displaystyle{B} & \\}$
\end{center}

%%%%%%%%%%%%%%%%%%%%%%%%%%%%%%%%%%%%%%%%%%%%%%%%%%%%%%%%%%%%%%%%%%%%%%%%%
%%%%%%%%%%%%%%%%%%%%%%%%%%%%%%%%%%%%%%%%%%%%%%%%%%%%%%%%%%%%%%%%%%%%%%%%

\end{definition}

%%%%%%%%%%%%%%%%%%%%%%%%%%%%%%%%%%%%%%%%%%%%%%%%%%%%%
%%%%%%%%%%%%%%%%%%%%%%%%%%%%%%%%%%%%%%%%%%%%%%%%%%%%

\proof (Theorem.\ref{HR})

Given a small enough $\epsilon >0$ and $\eta $ such that $0<\eta \ll \epsilon \ll 1$ consider the two fibrations, $f|:B_{\epsilon }(0)\cap f^{-1}(S_{\eta}^{1})\to S_{\eta}^{1}$ and $\displaystyle{\frac{f}{\|f \|}: S_{\epsilon }\setminus K_{\epsilon }\to S^{1}}$. By a diffeomorphism in $\mathbb{R}^{2}$ we can consider $\frac{1}{\eta}f|:B_{\epsilon }(0)\cap f^{-1}(S_{\eta}^{1})\to S^{1}$. For each $x\in B_{\epsilon }(0)\setminus V$ consider $p(s)$ the flow through x of vector fields
$\alpha $ in $B_{\epsilon }(0)\setminus V$  given by Lemma.\ref{L2}. That is, $\displaystyle{p'(s)=\alpha(p(s))}$ with $p(s_{0})=x$. By item a) this flow is transverse to all small sphere of radium $0<\epsilon_{1} \leq \epsilon $; by item b) it is transverse for all Milnor tube $B_{\epsilon }(0)\cap f^{-1}(S_{\eta}^{1})$, for all $\eta >0$, sufficiently small, and finally the item c) says that if the flow of vector fields $\alpha $ starts in a fiber $\displaystyle{F_{c}=(\frac{f}{\|f\|})^{-1}(c)}$ it will stay there all time in which it is defined. Because, using Proposition.\ref{p1} and item (c) above, we have $\displaystyle{\frac{d}{ds}(\frac{f(\alpha(s))}{\|f(\alpha(s)\|)})}=0.$ It means that $\displaystyle{\frac{f(\alpha(s))}{\|f(\alpha(s)\|)}=c}$, for all s in some maximal interval.

Of course, that item (a) says that this flow shall intersect $S_{\epsilon }\setminus K_{\epsilon }$ in some point say $y$. Then, defining $h:B_{\epsilon}(0)\cap f^{-1}(S_{\eta}^{1}) \to S_{\epsilon }\setminus K_{\epsilon }$ by $h(x)=y$, this smooth diffeomorphism $h$ is what we are looking for. \endproof

\begin{corollary}

Let $f=(P,Q):\mathbb{R}^{m},0\to \mathbb{R}^{2},0$, $m\geq 2$, a polynomial quasi-homogeneous map germ, of same weight and degree. Suppose $f$ has isolated critical point at origem $0\in \mathbb{R}^{m}$. Then, the fibrations ${\bf (1)}$ and ${\bf (2)}$ are equivalents.

\end{corollary}

\proof In \cite{SM} the authors proved that a necessary and sufficient condition to existence of fibration {\bf (2)} is that the projection $\displaystyle{\frac{f}{\|f \|}: S_{\epsilon_{0}}\setminus K_{\epsilon_{0}}\to S^{1}}$ be a submersion for all $\epsilon >0$, small enough, which is equivalent to have $\displaystyle{P(x)\nabla Q(x)-Q(x)\nabla P(x)}$ and the vector position $\displaystyle{x}$ linearly independent in $\displaystyle{B_{\epsilon}(0)\setminus V}$. By Lemma.\ref{L1}, these two vectors are orthogonal in $\displaystyle{B_{\epsilon}(0)\setminus V}$ hence, the projection is submersion and the fibration ${\bf (2)}$ exist. Since we have isolated singularity at origin, by Milnor result \cite{Mi} we have the existence of fibration {\bf (1)}. Now, applying the main theorem we have that these two fibrations are equivalents. \endproof

\begin{remark}

A similar approach to prove that these two fibrations are equivalents appear in \cite{SM}, but in a slight different way.

\end{remark}

%%%%%%%%%%%%%%%%%%%%%%%%%%%%%%%%%%%%%%%%%%%%%%%%%%%
%%%%%%%%%%%%%%%%%%%%%%%%%%%%%%%%%%%%%%%%%%%%%%%%%%%%

\end{document}